\documentclass[11pt,a4paper,fleqn]{article}
\usepackage{amsthm,amsmath,amssymb,amsfonts}
\newtheorem{theorem}{Theorem}

\numberwithin{equation}{section}
\numberwithin{lemma}{section}
\numberwithin{theorem}{section}
\numberwithin{corollary}{section}

\allowdisplaybreaks

\begin{document}
	\title{On the matrix $q$-Kummer equation and its solutions}
	\author{Ravi Dwivedi\footnote{E-mail: dwivedir999@gmail.com}  
		\\  Govt. Naveen College Bhairamgarh\\ Bijapur (CG) 494450, India.}
	\maketitle
	\begin{abstract}
In the present paper, a general theory for the second-order matrix difference equation of bilateral type is discussed. We introduced the matrix $q$-Kummer equation of bilateral type and presented the $q$-Kummer matrix function as a series solution. Later, we obtain the integral solutions of this Kummer equation and solution at $\infty$. We also give a brief idea about the matrix Gauss difference equation of bilateral type and its solutions. 
 
		\medskip
		\noindent\textbf{Keywords}: Basic hypergeometric function, Difference equations, Matrix functional calculus.
		
		\medskip
		\noindent\textbf{AMS Subject Classification}: 15A15, 39A13, 33D15.
	\end{abstract}
	
	\section{Introduction}\label{s:i}
	Most of the special functions are presented as a solution of either differential equation or integral equation and hence they play a vital role in the study of special functions. Due to a wide range of applications the Gauss hypergeometric equation and Kummer hypergeometric equation obtained special treatment from mathematician and physicist. These equations are studied from various point of view, for example the $q$-analogue of these two equations are given in \cite{gr} and \cite{jcf} respectively; The matrix version of Gauss hypergeometric function can be seen in \cite{jjc98b, jc00} and its $q$-analogue in \cite{as14}.
	
	In this paper, alongside the theory of the bilateral type matrix differential equation, we have tried to develop the same for the bilateral type matrix difference equation. In particular, we give a general form of a second order matrix difference equation of bilateral type and analyse its solutions. Furthermore, we define the ordinary and singular points of the $q$-Kummer matrix equation of  bilateral type and then present solutions on its singular points. The treatment of the paper is as follows:
	
	In section~2, we list some basic definitions and preliminaries that are needed in the subsequent sections. In section~3, a general second order bilateral type matrix $q$-difference equation is discussed. In section~4, the $q$-Kummer matrix function is defined and its regions of convergence is deduced. We also give matrix $q$-Kummer equation and various form of solutions is discussed. At last, in section~5, we added a conclusion part which includes the scope of this initial work. 
	
	\section{Basic definitions and Preliminaries}
The generalized hypergeometric function is defined by \cite{emo,edr, slater1966} 
\begin{align}
{}_rF_s (z) = 	{}_rF_s (\lambda_1, \dots, \lambda_r; \mu_1, \dots, \mu_s; z) = \sum_{n\ge 0} \frac{(\lambda_1)_n \cdots (\lambda_r)_n}{(\mu_1)_n \cdots (\mu_s)_n} \, \frac{z^n}{n !},\label{1.1}
\end{align}
where 	$\lambda_1, \dots, \lambda_r, \mu_1, \dots, \mu_s, z \in \mathbb{C}$ such that $\mu_j\neq 0, -1, -2,\dots,  1 \le j \le s$ and the Pochhammer symbol $(\lambda)_n$ is defined as
\begin{equation}
(\lambda)_n = \begin{cases}
	1, & \text{ if } n = 0,\\
	\lambda (\lambda + 1) \cdots (\lambda + n - 1), & \text{ if } n \ge 1.
\end{cases}
\end{equation} 
 The infinite series \eqref{1.1} converges absolutely in the following regions
 \begin{enumerate}
 	\item For $|z|<\infty$ if $r\leq s$;
 	
 	\item For $|z|<1$ if $r = s + 1$;
 	
 	\item For $|z|= 1$ if $r = s + 1$ and $\left(\sum_{j=1}^{s} \mu_j - \sum_{j=1}^{r} \lambda_j\right) > 0$.
 \end{enumerate}
The differential equation satisfied by generalized hypergeometric function \eqref{1.1} is 
\begin{align}
	\left[\theta (\theta + \mu_1 - 1) \cdots (\theta + \mu_s - 1) - z (\theta + \lambda_1) \cdots (\theta + \lambda_r)\right] \, {}_rF_s (z) = 0,\label{1.3}
\end{align}
where $\theta = z \frac{d}{dz}$. 

A function $g_q$ with parameter $q$ is called the $q$-analogue of an ordinary function $g$ if $g_q$ approaches to $g$ whenever $q \rightarrow 1$. With this basic idea, one of the $q$-analogue of differential equation \eqref{1.3} has been given by Jia \emph{et. al.}, \cite{jcf}, as
\begin{align}
	[\delta_q & (q^{\mu_1 - 1} \delta_q + [\mu_1 - 1]_q) \cdots  (q^{\mu_s - 1} \delta_q + [\mu_s - 1]_q)\nonumber\\
	& - z \, (q^{\lambda_1} \delta_q + [\lambda_1]_q) \cdots (q^{\lambda_r} \delta_q + [\lambda_r]_q)] u (z) = 0,\label{1.4}
\end{align}
where $\delta_q u(z) = z D_q u(z) = z \, \frac{u(z) - u(qz)}{z - qz}$, $ [\lambda]_q = \frac{1 - q^\lambda}{1 - q}$. As $\lim_{q \rightarrow 1} \delta_q \rightarrow \theta$ and $\lim_{q \rightarrow 1} [\lambda]_q \rightarrow \lambda$, so the $q$-difference equation \eqref{1.4} will tend to the ordinary differential equation \eqref{1.3} when $q \rightarrow 1$. One can verify that the equation \eqref{1.4} has a solution
\begin{align}
	u(z) & = {}_r \phi_s (q^{\lambda_1}, \dots, q^{\lambda_r}; q^{\mu_1}, \dots, q^{\mu_s}; q; (1 - q)^{1 + s - r} z)\nonumber\\
	& = \sum_{n\ge 0} \frac{(q^{\lambda_1}; q)_n \cdots (q^{\lambda_r}; q)_n }{(q^{\mu_1}; q)_n \cdots (q^{\mu_s}; q)_n} \, \frac{(1 - q)^{n \, (1 + s - r)} z^n}{(q; q)_n},\label{1.5}
\end{align} 
where 
\begin{equation*}
	(q^{\lambda_1};q)_n= \begin{cases}
		1, & \text{if $n=0$,} \\
		(1-q^{\lambda_1})(1-q^{\lambda_1 + 1}) \cdots (1-q^{\lambda_1 + n - 1}), &\text{if $n\geq1$}.
	\end{cases}
\end{equation*}
The series in \eqref{1.5} converges absolutely for $\left\vert (1 - q)^{1 + s - r} z\right\vert < 1$. 

For $r = 2$, $s = 1$, $\lambda_1 = \lambda$, $\lambda_2 = \nu$ and $\mu_1 = \mu$, the equation \eqref{1.4} produce the $q$-analogue of Gauss hypergeometric equation as
\begin{align}
	[\delta_q \, (q^{\mu - 1} \delta_q + [\mu - 1]_q)  - z \, (q^{\lambda} \delta_q + [\lambda]_q) \, (q^{\nu} \delta_q + [\nu]_q)] \, u (z) = 0.
\end{align}
Further, for $r = 1$, $s = 1$, $\lambda_1 = \lambda$ and $\mu_1 = \mu$, the equation \eqref{1.4} leads to the $q$-analogue of Kummer equation
\begin{align}
[\delta_q \, (q^{\mu - 1} \delta_q + [\mu - 1]_q)  - z \, (q^{\lambda} \delta_q + [\lambda]_q)] \, u (z) = 0,
\end{align}
which can also be written as
\begin{align}
	q^{\nu} \, z \, D_q^2 \, u(z) + \left([\nu]_q - q^\lambda \, z\right) D_q \, u(z) - [\lambda]_q u(z) = 0. 
\end{align}
The detailed study on this particular $q$-Kummer equation has been done in \cite{jcf}.

For a complex valued function $f(z$), the $q$-derivative is defined by \cite{kac2002}
\begin{align}
	D_q f(z)  = \frac{f(z) - f(qz)}{(1 - q) z}
\end{align}
and the $q$-derivative of product of two complex valued functions $f(z)$ and $g(z)$ is given as
\begin{align}
	D_q (f(z).g(z)) & =  f(qz) \, D_q \, g(z)  + g(z) \, D_q \, f(z).\label{6.32}
\end{align}
By symmetry, we can interchange $f$ and $g$
\begin{align}
	D_q (f(z).g(z)) = f(z) \, D_q \, g(z)  + g(qz) \, D_q \, f(z).\label{6.33}
\end{align}
The $q$-derivative of quotient of $f(z)$ and $g(z)$ is given by
\begin{align}
	D_q \left(\frac{f(z)}{g(z)}\right) = \frac{g(z) \, D_q \, f(z) - f(z) \, D_q \, g(z)}{g (z) \, g(qz)}.
\end{align}
If we use \eqref{6.33} instead of \eqref{6.32}, we get
\begin{align}
	D_q \left(\frac{f(z)}{g(z)}\right) = \frac{g(qz) \, D_q \, f(z) - f(qz) \, D_q \, g(z)}{g (z) \, g(qz)}.
\end{align}
The definite $q$-integral of a function $g(u)$, for complex number $u$ and $z$, is defined by \cite{gr, kac2002}
\begin{align}
	\int_{0}^{z} g(u) \, d_q u = (1 - q) \, \sum_{j \ge 0} z \, q^j \, g(z q^j).
\end{align}
In general,
\begin{align}
	\int_{0}^{z} g(u) \, d_q \, h(u) = 	\int_{0}^{z} g(u) \, D_q h(u) \, d_q u = \sum_{j \ge 0} g (zq^j) \left(h(zq^j) - h(z q^{j + 1})\right). 
\end{align}
The $q$-analogue of by part integration is given by, \cite{gr,kac2002}
\begin{align}
	\int_{0}^{z} g(qu) \, d_q \, h(u) =  g(z) \, h(z)  - g(0) \, h(0)  -	\int_{0}^{z} h(u) \, D_q g(u). 
\end{align}
Further, the $q$-integral of $g(u)$ on $[0, \infty)$ is given by
\begin{align}
	\int_{0}^{\infty} g(u) \, d_q u = (1 - q) \, \sum_{j = -\infty}^{\infty} q^j \, g(q^j). 
\end{align}
The $q$-derivative of a definite $q$-integral is given by the following formula, \cite{jcf}
\begin{align}
	D_q \left(\int_{0}^{\frac{1}{a z^k}} g(z, u) \, d_q u \right) =  \left(\int_{0}^{\frac{1}{a z^k}} D_q g(z, u) \, d_q u \right) - \sum_{j = 0}^{k - 1} \frac{q^{j - k} }{a \, z^{k + 1}} g \left(qz, \frac{q^{j - k}} {a z^k}\right),\label{2.21}
\end{align}
where $a$, $z$ are complex numbers and $D_q$ is the $q$-derivative with respect to $z$. Letting $q \rightarrow 1$, the formula \eqref{2.21} yields
\begin{align}
	\frac{d}{dz} \left(\int_{0}^{\frac{1}{a z^k}} g(z, u) \, d_q u \right) = \int_{0}^{\frac{1}{a z^k}} \frac{\partial}{\partial z} g(z, u) \, d_q u - \frac{k}{a \, z^k} \, g(z, \frac{1}{a \, z^k}). 
\end{align}
Let $P$ be a $p$-square complex matrix. Then,  $\alpha(P) = \max\{\, \Re(z) \mid z\in \sigma(P)\, \}$ and $\beta(P) = \min\{\, \Re(z) \mid z \in \sigma(P)\, \}$,  where $\Re(z)$  is the real  part of $z \in \mathbb{C}$ and $\sigma(P)$ is the set of all eigenvalues of $P$.  The matrix $P \in \mathbb{C}^{p \times p}$ is called positive stable if  $\beta(P) >0$. 

The $q$-analogue of a matrix $P$ is defined by \cite{as12}
\begin{align}
	[P]_q = \frac{I - q^P}{1 - q}, \quad q^P = \exp (P \, \ln q), 
\end{align}  
where $\ln$ is the principal branch of logarithmic function. The matrix analogue of $q$-shifted factorial function, for $P \in \mathbb{C}^{p \times p}$, is given by \cite{as12, as14}
\begin{equation*}
	(q^{P};q)_n= \begin{cases}
		I, & \text{if $n=0$,} \\
		(I-q^{P})(I-q^{P + I}) \cdots (I-q^{P + (n - 1)I}), &\text{if $n\geq1$}.
	\end{cases}
\end{equation*}
Furthermore, we have 
\[(q^P;q)_\infty = \lim_{n \rightarrow \infty} \prod_{k=0}^{n-1} (I-q^{P + kI}) = \prod_{k=0}^{\infty} (I-q^{P + kI}),\]
which is convergent for $\vert q\vert < 1$, \cite{as12}. If $\Vert q^P\Vert < 1$ and $q^{-k} \notin \sigma(q^{P} )$, $k \in \mathbb{Z}^+$ (set of positive integers), then $(q^P;q)_\infty$ converges invertibly, \cite{as12, as14, wt95}.

The $q$-gamma matrix function, for $P\in \mathbb{C}^{p\times p}$ such that $\beta (P) > 0$, is defined as \cite{as12} 
\begin{equation}
	\Gamma_q(P) = \int_{0}^{\frac{1}{1-q}} u^{P-I} E_q(-qu) \ d_qu,
\end{equation}
where $E_q(z)$ is the $q$-exponential function defined by
\begin{align}
E_q(z) = 	\sum_{j = 0}^{\infty}  \frac{q^{j \choose 2} \, \left(z\right)^j}{[j]_q !} = (- (1 - q)  z ; q)_\infty.
\end{align}
Furthermore, if $P$ and $Q$ commute, then the $q$-beta matrix function is given by \cite{as12}
\begin{align}
	\mathcal{B}_q(P, Q) = \mathcal{B}_q(Q, P) & = \int_{0}^{1} (uq ; q)_\infty \ (uq^Q ; q)_\infty ^{-1} \ u^{P-I} d_qu \nonumber\\
	&  = \Gamma_q (P) \ \Gamma_q (Q) \ \Gamma _q ^{-1} (P + Q), 
\end{align}
where inverse of $q$-gamma matrix function $\Gamma_q^{-1} (P)$ is defined by
\begin{equation}
	\Gamma _q ^{-1} (P) = [P]_q \, [P+I]_q \, \cdots \, [P+(n-1)I]_q \, \Gamma _q ^{-1} (P+nI), \quad n = 1,2, \dots. \label{eq14} 
\end{equation}
In view of the equation \eqref{eq14}, one can write
\begin{equation}
	(q^P ; q)_n = (1-q)^n \, \Gamma_q^{-1}(P) \, \Gamma _q (P+nI), \quad n = 1,2, \dots. \label{eq15} 
\end{equation}
For detailed study on basic matrix functions, one can check these papers and references therein \cite{ds2, ds8, as12, as14}.

\section{Matrix $q$-difference equations}
Consider the following general matrix $q$-difference equation 
\begin{align}
	[\delta_q & (q^{Q_1 - 1} \delta_q + [Q_1 - 1]_q) \cdots  (q^{Q_s - 1} \delta_q + [Q_s - 1]_q)\nonumber\\
	& - z \, (q^{P_1} \delta_q + [P_1]_q) \cdots (q^{P_r} \delta_q + [P_r]_q)] W (z) = 0,\label{1.9}
\end{align}
where $P_1, \dots, P_r, Q_1, \dots, Q_s$ are matrices in $\mathbb{C}^{p \times p}$ such that $q^{-k} \notin \sigma(q^{Q_j} ), 1 \le j \le s$, $k \in \mathbb{Z}^+$. One can check that the generalized basic hypergeometric matrix function, defined by 
\begin{align}
 & {}_r \phi_s (q^{P_1}, \dots, q^{P_r}; q^{Q_1}, \dots, q^{Q_s}; q; (1 - q)^{1 + s - r} z)\nonumber\\
	& = \sum_{n\ge 0} (q^{P_1}; q)_n \cdots (q^{P_r}; q)_n (q^{Q_1}; q)^{-1}_n \cdots (q^{Q_s}; q)^{-1}_n \, \frac{(1 - q)^{n \, (1 + s - r)} z^n}{(q; q)_n},\label{1.11}
\end{align} 
is the solution of the matrix $q$-difference equation \eqref{1.9}. While proving this, we have to  directly or indirectly assume that the matrices $P_1, \dots, P_r, Q_1, \dots, Q_s$ are commutative. To reduce this comprehensive assumption  of commutativity, we introduce the concept of matrix $q$-difference equations of bilateral type. 

The  bilateral type matrix differential equation has been introduced by J\'odar and Cort\'es \cite{jjc98b, jc00} and studied the Gauss hypergeometric matrix equation systemically. Author extended this work for several variable hypergeometric matrix functions \emph{viz.} Appell matrix functions, Lauricella matrix functions of three and $n$ variables \cite{ds1, ds4, ds5}. This motivates to study second order matrix $q$-difference equation of bilateral type in a systematic way. Consider the second order matrix $q$-difference equation of bilateral type
\begin{align}
	&\phi_{1, q} (z) \, D_q^2 U(z) + D_q^2 U(z) \, \phi_{2, q} (z) + \phi_{3, q} (z) \, D_q U(z) + D_q U(z) \, \phi_{4, q} (z)\nonumber\\
	& \quad  + \phi_{5, q} (z) \,D_q U(z) \, \phi_{6, q} (z) + \phi_{7, q} (z) \,  U(z) \, \phi_{8, q} (z) = 0,\label{1.12}
\end{align}
where $\phi_{i, q} (z), 1 \le i \le 8$ are bounded continuous functions in an open disc $B(z_0, p)$ of radius $p$ centred  at $z_0$ given by $\phi_{i, q} : B(z_0, p) \rightarrow E_q$. Here $E_q$ is the $q$-Banach space with $2$-norm such that $\lim_{q \rightarrow 1} E_q = E$ (Banach space of all $p$-square matrices). 

The point $z = z_0$ is called an ordinary point of $q$-difference  equation \eqref{1.12} if $\phi_{1, q} (z_0)$ and $\phi_{2, q} (z_0)$ are invertible. In case, the inverse of  $\phi_{1, q} (z_0)$ and $\phi_{2, q} (z_0)$ do not exist but each $\phi_{j, q} (z_0), 3 \le j \le 8$ is invertible, then we called $z = z_0$ a singular point. Further, if limits
\begin{align}
	\lim_{z \rightarrow z_0} (z - z_0) \, \phi_{j, q} (z) \, [\phi_{i, q} (z)]^{-1}, \quad i = 1, 2, \ j = 3, 4, 5, 6
\end{align} 
and 
\begin{align}
	\lim_{z \rightarrow z_0} (z - z_0)^2 \, \phi_{k, q} (z) \, [\phi_{i, q} (z)]^{-1}, \quad i = 1, 2, \ k = 7, 8
\end{align} 
exist, then the point $z = z_0$ is called a regular singular point of \eqref{1.12}. If any of the limits given above do not exist, then the singular point is irregular.  

 Let $\phi_q: B (z_0; p) \times E_q \rightarrow E_q$ be a continuous function such that 
	\begin{align}
		\Vert \phi_q (z, X_{1, q}) - \phi_q (z, X_{2, q}) \Vert \le K \, (\vert z - z_0\vert) \, \Vert  X_{1, q} -  X_{2, q}\Vert,
	\end{align}
where $z \in B(z_0, p)$; $X_{1, q}, X_{2, q} \in E_q$  and $K_q$  is a real valued continuous function on $[0, p]$. Then,  there exists an unique solution $W_q(z)$ of 
\begin{equation}
	D_q (X_q) = \phi_q (z, X_q),
\end{equation}
 such that $W_q (z_0) = X_{0, q}$. 
 
 Any two solutions $W_{1, q}$ and $W_{2, q}$ of second order $q$-difference equation \eqref{1.12} are called fundamental set of solution if there exist a solution $W_q$ of \eqref{1.12} which has a unique representation in the form
 \begin{align}
 	W_q (z) = W_{1, q}(z) \, A + W_{2, q} (z) \, B, 
 \end{align}
 where matrices $A$ and $B$ are uniquely obtained by $W_q(z)$.

\section{The $q$-Kummer matrix equation}
Before discussing $q$-Kummer matrix equation and its solution we will define $q$-Kummer matrix function and give its regions of convergence. Let $S$, $T$ be $p$-square matrices such that $q^{-k} \notin \sigma(q^T), k \in \mathbb{Z}^+$. Then, we define the $q$-Kummer matrix function as
\begin{align}
	{}_1\phi_1(q^S; q^T; q; z) = \sum_{n\ge 0} \frac{(q^S; q)_n (q^T; q)_n^{-1}}{(q; q)_n} \, z^n.
\end{align}
For convergence,
\begin{align}
\Vert {}_1\phi_1(q^S; q^T; q; z)\Vert & \le  \sum_{n\ge 0} \left \Vert \frac{(q^S; q)_n (q^T; q)_n^{-1}}{(q; q)_n} \, z^n\right \Vert\nonumber\\
& \le \sum_{n\ge 0} \frac{(-\Vert q^S\Vert; \vert q\vert)_n \Vert (q^T; q)_n^{-1}\Vert}{\vert (q; q)_n \vert} \vert z\vert^n = \sum_{n\ge 0} U_n \, \vert z\vert^n.
\end{align}
Now, using the ratio test
\begin{align}
\lim_{n \rightarrow \infty} \left \vert \frac{U_{n + 1}}{U_n}\right\vert	& =  \lim_{n \rightarrow \infty} \left \vert \frac{(-\Vert q^S\Vert; \vert q\vert)_{n+1} \Vert (q^T; q)_{n + 1}^{-1}\Vert \vert z\vert^{n + 1} \vert (q; q)_{n} \vert}{(-\Vert q^S\Vert; \vert q\vert)_{n} \Vert (q^T; q)_{n}^{-1}\Vert \vert z\vert^{n} \vert (q; q)_{n+ 1} \vert}\right\vert\nonumber\\
	& = \lim_{n \rightarrow \infty} \left \vert \frac{(1 + \vert q\vert^n \Vert q^S\Vert) \Vert (I - q^n q^T)^{-1} \Vert \vert z\vert}{1 - \vert q\vert ^{n + 1}}\right\vert\nonumber\\
	& \le \lim_{n \rightarrow \infty} \frac{(1 + \vert q\vert^n \Vert q^S\Vert) \, \vert z\vert}{(1 - \vert q\vert ^{n + 1}) \, (1 - \vert q\vert^n \Vert q^T\Vert)}  = \vert z \vert.
\end{align}
Thus the $q$-Kummer matrix function ${}_1\phi_1(q^S; q^T; q; z)$ converges absolutely when $\vert q\vert < 1$ and $\vert z\vert < 1$. 

Next, we give the matrix $q$-Kummer equation of bilateral type and discuss the distinct solutions.  Let $S$, $T$ be  matrices in $\mathbb{C}^{p \times p}$ such that $q^{-k} \notin \sigma(q^T), k \in \mathbb{Z}^+$. Then, we introduce the matrix $q$-Kummer equation of bilateral type as 
\begin{align}
	z D_q^2 U(z) q^T + D_q U(z) [T]_q - q^S z D_q U(z) - [S]_q U(z) = 0.\label{2.4}
\end{align}
This is an equation of the form \eqref{1.12}, on comparing we have
\begin{align}
	\phi_{2, q} (z) = z \, q^T, &\quad \phi_{3, q} (z) = - z \, q^S, \quad 		\phi_{4, q} (z) = [T]_q, \quad \phi_{7, q} (z) = - [S]_q.
\end{align}
The only singular points of the equation \eqref{2.4} are $z=0$ and $z=\infty$. First, we give our attention in finding the series solution of the  matrix $q$-difference equation about the singular point $z = 0$. Then we look into for another form of solutions \emph{viz.} solution at $\infty$ and integral solutions. Let us assume the series solution of \eqref{2.4} be in the form 
\begin{align}
	U(z) = \sum_{n\ge 0} U_n z^n
\end{align}
which yields
\begin{align}
	D_q \, U(z) = \sum_{n\ge 1} U_n \, [n]_q \,  z^{n - 1}, \quad D_q^2 \, U(z) = \sum_{n\ge 2} U_n \, [n]_q \, [n - 1]_q \,  z^{n - 2}. 
\end{align}
Now, substituting these values in \eqref{2.4}, we get
\begin{align}
&\sum_{n\ge 2} U_n \, [n]_q \, [n - 1]_q \,  z^{n - 1} q^T + \sum_{n\ge 1} U_n \, [n]_q \,  z^{n - 1} [T]_q - q^S \, \sum_{n\ge 1} U_n \, [n]_q \,  z^{n}\nonumber\\
& \quad  - [S]_q \sum_{n\ge 0} U_n \,  z^{n} = 0.\label{2.7}
\end{align}
After some manipulation, the equation \eqref{2.7} gives 
\begin{align}
	& \sum_{n\ge 2} \left(U_{n + 1} [n + 1]_q \, [n]_q q^T + U_{n + 1} [n + 1]_q [T]_q - q^S U_n [n]_q - [S]_q U_n \right) \, z^n \nonumber\\
	& \quad  + (U_2 [2]_q q^T + U_2 [2]_q [T]_q - q^S U_1 - [S]_q U_1) \, z + U_1 [T]_q - [S]_q U_0 = 0. 
\end{align}
Equating the powers of $z$ and hence one can obtain
\begin{align}
	U_1 \, [T]_q = [S]_q U_0 \implies U_1 = \frac{[S]_q \, U_0 \, [T]_q^{-1} \, (1 - q)}{(1 - q)};
\end{align}
\begin{align}
	U_2  =  \frac{[S + I]_q \, U_1 \, [T + I]_q^{-1}}{[2]_q} & = \frac{[S + I]_q \, [S]_q \, U_0 \, [T + I]_q^{-1} \, [T]_q^{-1}}{[2]_q}\nonumber\\
	& =  \frac{[S + I]_q \, [S]_q \, U_0 \, [T + I]_q^{-1} \, [T]_q^{-1} \, (1 - q)^2}{(1 - q) \, (1 - q^2)}.
\end{align}
In general
\begin{align}
	U_{n + 1} & = \frac{[S + nI]_q \, U_n \, [T + nI]_q^{-1}}{[n + 1]_q}\nonumber\\
	& =  \frac{[S + nI]_q \, \cdots [S]_q \,  U_0 \, [T + nI]_q^{-1} \, \cdots [T]_q^{-1} \, (1 - q)^n}{(1 - q) \, (1 - q^2) \cdots (1 - q^{n + 1})}.
\end{align}
Hence there exists a solution $U_1 (z)$ of \eqref{2.4}, well defined in $\vert (1 - q) \, z\vert < 1$ and satisfying $U_1 (0) = I$, is given by 
\begin{equation}
	U_1(z) = {}_1\phi_1(q^S; q^T; q; (1 - q) \, z).\label{4.13}
\end{equation}
Since equation \eqref{2.4} is a second order matrix difference equation, so it must have another series solution.  Suppose, for $ST = TS$, the second solution has the form
\begin{align}
	U(z) = V(z) \, z^{I - T}, \quad z^{I - T} = \exp ((I - T) \ln z).\label{2.13}
\end{align} 
The $q$-derivatives of \eqref{2.13} are given by
\begin{align}
	D_q U(z) & = V(z) \, [I - T]_q z^{-T} + D_q V(z) (qz)^{I - T}\label{2.14}
\end{align}
\begin{align}
	D_q^2 U(z) & = V(z) \, [I - T]_q \, [-T]_q \, z^{-T - I} + D_q V(z) \, [I - T]_q \, (qz)^{-T} + D_q^2 V(z) (q^2 z)^{I - T}\nonumber\\
	& \quad  + q D_q V (z) \, [I - T]_q \, (qz)^{-T}.\label{2.15}
\end{align}
Substituting these values in \eqref{2.4} yield
\begin{align}
	z \, D_q^2 V(z) \, q^{2I - T} + D_q V(z) \, [2I - T]_q - z \, q^{S + I - T} \, D_q V(z) - [S + I - T]_q V(z) = 0.\label{2.16}
\end{align}
The equation \eqref{2.16} is a matrix $q$-difference equation of type \eqref{2.4} and its solution is given by ${}_1\phi_1(q^{\tilde{S}}; q^{\tilde{T}}; q; (1 - q) \, z)$, where $\tilde{S} = S + I - T$ and $\tilde{T} = 2I - T$. Thus we can conclude that the matrix $q$-Kummer equation \eqref{2.4} has two fundamental solutions as
\begin{align}
	U_1(z) & = {}_1\phi_1(q^S; q^T; q; (1 - q) \, z);\nonumber
	\\[5pt]
	U_2 (z) & = z^{I - T} \, {}_1\phi_1(q^{S + I - T}; q^{2I - T}; q; (1 - q) \, z), \quad ST = TS. 
\end{align} 
In view of the convergence condition for $q$-Kummer matrix function, the solutions $U_1(z)$ and $U_2(z)$ converges absolutely in the region $\vert q\vert < 1$ and $\vert (1 - q) \, z\vert < 1$.

\subsection{Solutions at $\infty$}
Now, we seek for the solution of matrix  $q$-Kummer equation at infinity and so we rewrite the equation \eqref{2.4} as
\begin{align}
	U(q^2 z) q^{T - I} - [I + q^S (q - 1) z] \, U(qz) - U(qz) \, q^{T - I} + [1 + (q - 1)z] \, U(z) = 0. \label{2.18}
\end{align}
Let $u = z^{-1}$, $s = q^{-1}$, $W(u) = U(z^{-1})$ and $ST = TS$. Then \eqref{2.18} becomes
\begin{align}
&	u \, W(s^2 u) - [(1 - s) \, s^{T - S - 2I} + (I + s^{T - I}) \, u] \, W(su)\nonumber\\
	& \quad  + [(1 - s) \, s^{T - 2I} + u \, s^{T - I}] \, W(u) = 0. 
\end{align}
Which can be written in the derivative form as
\begin{align}
	u^3 \, D_s^2 W(u) + [u\, s^{T - S - 3I} - u^2 \, [T - 2I]_s ] D_s \, W(u) + s^{T - 3I} \, [-S]_s \, W(u) = 0,\label{2.20}
\end{align}
where $D_s \, W(u) = \frac{W(u) - W(su)}{(1 - s) \, u}$ and $[S]_s = \frac{I - s^S}{1 - s}$. This is the matrix $s$-difference equation of the form \eqref{1.12} where
\begin{align}
	\phi_{1, s} (u) = u^3, \ \phi_{3, s} (u) = u\, s^{T - S - 3I} - u^2 \, [T - 2I]_s, \ \phi_{7, s} (u) = s^{T - 3I} \, [-S]_s.
\end{align}
All of the above are well defined near $u = 0$ except  the limit
\begin{align}
	\lim_{u \rightarrow 0} \frac{u \, \phi_{3, s} (u)}{\phi_{1, s} (u)} = 	\lim_{u \rightarrow 0} \frac{u \, [u\, s^{T - S - 3I} - u^2 \, [T - 2I]_s]}{u^3}. 
\end{align}
Thus, point $u = 0$ is irregular singular for \eqref{2.20} and hence the point $z = \infty$ is irregular for \eqref{2.4}. Therefore, the equation \eqref{2.4} does not have a solution at $\infty$. Unfortunately, we are not getting a convergent solution of \eqref{2.4} at $z = \infty$ but it is possible to get some integral solutions which are being discussed in the next section.  

\subsection{Integral solutions of $q$-Kummer matrix equation}
Consider one of the integral solution of \eqref{2.4} in the form
\begin{align}
	U_1(z) = \int_{0}^{\frac{1}{(1 - q) q^2 z}} E_q^{- q z u} \, f(q u) \, d_q u, 
\end{align}
where $E_q^{- q z u}$ and its $q$-derivative are defined by \cite{gr}
\begin{align}
E_q^{- q z u}  = ((1 - q) q z u; q)_\infty \text{ and }
D_q E_q^{- q z u}  = - qu \, E_q^{- q^2 z u}. \label{e2.34}
\end{align}  
Using \eqref{2.21} and \eqref{e2.34}, one can get 
\begin{align}
	&D_q \,   U_1 (z)\nonumber\\
	 & = \int_{0}^{\frac{1}{(1 - q) q^2 z}} (- qu) \, E_q^{- q^2 z u} \, f(qu) \, d_q u - \left(\frac{1}{(1 - q) q^3 z^2} \, E_q^{- q^2 z u} f(qu)\right)_{u = 1/ (1 - q)q^3 z}\nonumber\\
	& = - \int_{0}^{\frac{1}{(1 - q) q^2 z}} qu \, E_q^{- q^2 z u} \, f(qu) \, d_q u,
\end{align}
and 
\begin{align}
	D_q^2 \, U_1 (z) & = - D_q \left(\int_{0}^{\frac{1}{(1 - q) q^2 z}} qu \, E_q^{- q^2 z u} \, f(qu) \, d_q u\right)\nonumber\\
	& = \int_{0}^{\frac{1}{(1 - q) q^2 z}} q^3 u^2 \, E_q^{- q^3 z u} \, f(qu) \, d_q u.
\end{align}
Substituting these values in \eqref{2.4}, we get
\begin{align}
& \int_{0}^{\frac{1}{(1 - q) q^2 z}} [-f(qu) D_{q, u} (E_q^{- q^2 z u}) q^{T + I} u^2 - u \, q^S \, f(q u) \, D_{q, u} (E_q^{- q z u}) \nonumber\\
& \quad - q u f(q u) E_q^{- q^2 z u}  [T]_q - [S]_q f(qu) E_q^{- q z u}] d_q u\nonumber\\
& = \int_{0}^{\frac{1}{(1 - q) q^2 z}}  E_q^{- q^2 z u} \, D_{q, u} (f(u) q^{T- I} u^2) + E_q^{- q z u} D_{q, u} (q^{S - I} f(u) \, u) \nonumber\\
& \quad - q u f(q u) E_q^{- q^2 z u}  [T]_q - [S]_q f(qu) E_q^{- q z u}] d_q u\nonumber\\
& = \int_{0}^{\frac{1}{(1 - q) q^2 z}} E_q^{- q^2 z u} \{D_{q, u} (f(u) q^{T- I} u^2) + (1 - (1 - q) q z u) D_{q, u} (q^{S - I} f(u) \, u) \nonumber\\
& \quad - q u f(q u)   [T]_q - [S]_q f(qu)  + (1 - q) q z u [S]_q f (q u)  \} d_q u = 0,
 \end{align}
where $D_{q, u}$ is the $q$-derivative with respect to $u$. Since the integral value is zero so one of the possibility is that
\begin{align}
&	D_{q, u} (f(u) q^{T- I} u^2) + (1 - (1 - q) q z u) D_{q, u} (q^{S - I} f(u) \, u) - q u f(q u)   [T]_q\nonumber\\
	&\quad   - [S]_q f(qu)  + (1 - q) q z u [S]_q f (q u)  = 0.
\end{align} 
Expanding the term containing $q$-derivative, we get a difference equation in terms of $f(u)$ as
\begin{align}
	\{1 + [q - (1 - q) qz] \, u\} f(qu) -q^{S - I} \{1 - (1 - q) qzu\} f(u) - u f(u) q^{T - I} = 0.\label{2.28} 
\end{align}
One can easily verify that for $ST = TS$ the matrix function, 
\begin{align}
	F(u) = u^{S - I} \, ([(1 - q) qz - q] u; q)_\infty \, ([(1 - q) qz - q^{T - S}] u; q)_\infty^{-1},
\end{align}
is a solution of recurrence relation \eqref{2.28}. Thus $U_1 (z)$ can be expressed as
\begin{align}
	U_1(z) & = \int_{0}^{\frac{1}{(1 - q) q^2 z}} E_q^{- q z u} \, {(qu)}^{S - I}  ([(1 - q) qz - q] q u; q)_\infty \nonumber\\
	& \quad \times   ([(1 - q) qz - q^{T - S}] q u; q)_\infty^{-1} \, d_q u. 
\end{align}
This is an integral solution of $q$-Kummer matrix equation \eqref{2.4}. We already discussed,  a second order difference equation must have two fundamental solutions. So there is an opportunity of getting another integral solution. To get the second solution we use the another $q$-analogue of exponential function exists in the literature \cite{gr}. We will look the second solution in the form
\begin{align}
	U_2 (z) = \int_{0}^{\infty} e_q^{-zu} \, f(q u) \, d_q u,
\end{align}
where 
\begin{align}
	e_q^{-z u} = \frac{1}{(-(1 - q) z u; q)_\infty}.
\end{align}
The $q$-derivative of $e_q^{-z u}$ with respect to $z$ is 
\begin{align}
 D_q e_q^{-z u} = - u \, e_q^{-z u}.\label{2.34}
\end{align}
Further, using \eqref{2.21} and \eqref{2.34}, one can get
\begin{align}
	D_q U_2 (z) = - \int_{0}^{\infty} u \, e_q^{-zu} \, f(q u) \, d_q u,
	\\[5pt]
	D_q^2 U_2 (z) =  \int_{0}^{\infty} u^2 \, e_q^{-zu} \, f(q u) \, d_q u.
\end{align}
Substituting these values in \eqref{2.4}, we have
\begin{align}
	\int_{0}^{\infty} f(q u) [- (D_{q, u} e_q^{-zu}) q^T u^2 - (D_{q, u} e_q^{-zu}) q^S u - ([S]_q + [T]_q u) e_q^{-zu}] d_q u = 0.
\end{align}
Which can be rewritten as
\begin{align}
\int_{0}^{\infty} e_q^{-zu} [D_{q, u} (f(u) q^{T - 2I} u^2) + D_{q, u} (f(u) q^{S - I} u) - f(q u) ([S]_q   + [T]_q u)] d_q u = 0.	
\end{align}
From above equation, one can obtain a matrix $q$-difference equation in terms of $f(u)$ as
\begin{align}
	(1 + u) f (q u) - (q^{S - I} + q^{T - 2I} u) f(u) = 0.\label{2.38}
\end{align} 
It is easy to verify that the matrix $q$-difference equation \eqref{2.38} has solution 
\begin{align}
f(u) = u^{S - I} {(-u; q)_\infty} (-q^{S - T - I} u; q)_\infty^{-1}.
\end{align}
Thus, $U_2 (z)$ can be expressed as
\begin{align}
	U_2 (z) = q^{S - I}  \, \int_{0}^{\infty} u^{S - I} \,  e_q^{-zu} \, {(-qu; q)_\infty} (-q^{S - T } u; q)_\infty^{-1} \, d_q u.
\end{align} 
We now summarize the integral solutions of \eqref{2.4} in the following theorem.
\begin{theorem}
	Let $S$ and $T$ be commuting  matrices in $\mathbb{C}^{p \times p}$ such that $S$, $T$ and $T - S$ are positive stable. Then the matrix equation \eqref{2.4} has two integral solutions as
	\begin{align}
		U_1(z) & = \int_{0}^{\frac{1}{(1 - q) q^2 z}} E_q^{- q z u} \, {(qu)}^{S - I}  ([(1 - q) qz - q] q u; q)_\infty \nonumber\\
	& \quad \times   ([(1 - q) qz - q^{T - S}] q u; q)_\infty^{-1} \, d_q u;\nonumber
	\\[5pt]
		U_2 (z) &= q^{S - I}  \, \int_{0}^{\infty} u^{S - I} \,  e_q^{-zu} \, {(-qu; q)_\infty} (-q^{S - T } u; q)_\infty^{-1} \, d_q u.
	\end{align}
\end{theorem}

\section{Conclusion} 
In this paper, we determined the  solutions of $q$-Kummer matrix equation of bilateral type. We discussed the series solutions, solution at infinity and integral solutions. One can observe that the only solution for non-commutative $S$ and $T$ is in \eqref{4.13}. For other solutions, we are directly or indirectly assuming that $S$ and $T$ are commuting. 

The present work of  ${}_1\phi_1 (q^S; q^T; q; (1 - q)z)$ can be easily extended for several variable basic hypergeometric matrix function. One can introduce the matrix $q$-difference equation of bilateral type and obtain its distinct form of solutions. In particular, let $P$, $Q$, $R$ be matrices in $\mathbb{C}^{p \times p}$ such that $q^{-k} \notin \sigma(q^R), k \in \mathbb{Z}^+$. Then, it is easy to show that the matrix $q$-difference equation,
\begin{align}
	\left[\delta_q (q^{R - I} \, \delta_q + [R - I]_q) - z (q^P \, \delta_q + [P]_q) \, (q^Q \, \delta_q + [Q]_q)\right] W(z) = 0,\label{3.1}
\end{align}
has basic Gauss hypergeometric matrix function ${}_2\phi_1 (q^P, q^Q; q^R; q; z)$ as a solution provided the matrices $P$, $Q$, $R$ commutes with each other. To overcome with such comprehensive assumption of commutativity, one can introduce a second order bilateral type matrix $q$-difference equation 
\begin{align}
	& z D_q^2 \, W(z) \, q^R - z [P]_q \, D_q W(z) \, q^Q + D_q W(z) \, [R]_q - z \, q^P \, D_q W(z) [Q + I]_q\nonumber\\
	& \quad  - z^2 q^{P + I} D_q^2 W(z) \,  q^Q - [P]_q \, W(z) \, [Q]_q = 0.   \label{3.2}
\end{align}
The equation \eqref{3.2} is a $q$-analogue of matrix Gauss hypergeometric equation given in \cite{jjc98b}. On comparing the above equation with \eqref{1.12}, we get
\begin{align}
	\phi_{1, q} (z) & = - z^2 \, q^{P + I}, \ \phi_{2, q} (z) = z \, q^{R - Q}, \ \phi_{3, q} (z) = - z \, [P]_q, \ \phi_{4, q} (z) = [R]_q \, q^{-Q}, \nonumber\\
	\phi_{5, q} (z) & = - z \, q^P, \ \phi_{6, q} (z) = [Q + I]_q \, q^{-Q}, \ \phi_{7, q} (z) = -[P]_q, \ \phi_{8, q} (z) = [Q]_q \, q^{-Q}.
\end{align} 
Treating in the same way to matrix $q$-Kummer equation, one can obtain that the  series solutions of the equation \eqref{3.2} are
\begin{align}
	W_1 (z) & = {}_2\phi_1 (q^P, q^Q; q^R; q; z), \quad QR = RQ;\nonumber
	\\[5pt]
	W_2 (z) & = z^{I - R} \, {}_2\phi_1 (q^{P + I - R}, q^{Q + I - R}; q^{2I - R}; q; z), \quad PR = RP, QR = RQ.
\end{align}
Finding the other form of solutions of matrix difference equation \eqref{3.2} is an interesting problem and yet to be done. To find the matrix $q$-difference equation satisfied by basic Gauss hypergeometric matrix function for non-commutative matrices is also an open problem and may be obtained in my upcoming articles.   
%
%

\end{document}